\documentclass[english,12pt]{smfart}
\usepackage{amssymb}
\usepackage{amsfonts}
\usepackage{amscd}
\usepackage[all]{xypic}

\swapnumbers

\newtheorem{theorem}[subsection]{Theorem}
\newtheorem{lem}[subsection]{Lemma}

\newtheorem{prop}[subsection]{Proposition}

\theoremstyle{definition}

\newtheorem{def-prop}[subsubsection]{Definition-Proposition}

\theoremstyle{remark}
\newtheorem{remark}[subsection]{Remark}

\theoremstyle{plain}

\numberwithin{equation}{subsection}

\def\boxit#1#2{\setbox1=\hbox{\kern#1{#2}\kern#1}%
\dimen1=\ht1 \advance\dimen1 by #1
\dimen2=\dp1 \advance\dimen2 by #1
\setbox1=\hbox{\vrule height\dimen1 depth\dimen2\box1\vrule}%
\setbox1=\vbox{\hrule\box1\hrule}%
\advance\dimen1 by .4pt \ht1=\dimen1
\advance\dimen2 by .4pt \dp1=\dimen2 \box1\relax}

\renewcommand{\theequation}{\thesubsection.\arabic{equation}}

\let\cal\mathcal

\def\AA{{\mathbf A}}

\def\CC{{\mathbf C}}

\def\FF{{\mathbf F}}
\def\GG{{\mathbf G}}

\def\NN{{\mathbf N}}

\def\PP{{\mathbf P}}
\def\QQ{{\mathbf Q}}
\def\RR{{\mathbf R}}

\def\cA{{\mathcal A}}

\def\cE{{\mathcal E}}

\def\cG{{\mathcal G}}
\def\cH{{\mathcal H}}

\def\cL{{\mathcal L}}

\def\cS{{\mathcal S}}

\def\cW{{\mathcal W}}

\mathchardef\alphag="7C0B
\mathchardef\betag="7C0C
\mathchardef\gammag="7C0D
\mathchardef\deltag="7C0E
\mathchardef\varepsilong="7C22
\mathchardef\varphig="7C27
\mathchardef\psig="7C20
\mathchardef\zetag="7C10
\mathchardef\epsilong="7C0F
\mathchardef\rhog="7C1A
\mathchardef\taug="7C1C
\mathchardef\upsilong="7C1D
\mathchardef\iotag="7C13
\mathchardef\thetag="7C12
\mathchardef\pig="7C19
\mathchardef\sigmag="7C1B
\mathchardef\etag="7C11
\mathchardef\omegag="7C21
\mathchardef\kappag="7C14
\mathchardef\lambdag="7C15
\mathchardef\mug="7C16
\mathchardef\xig="7C18
\mathchardef\chig="7C1F
\mathchardef\nug="7C17
\mathchardef\varthetag="7C23
\mathchardef\varpig="7C24
\mathchardef\varrhog="7C25
\mathchardef\varsigmag="7C26
\mathchardef\Omegag="7C0A
\mathchardef\Thetag="7C02
\mathchardef\Sigmag="7C06
\mathchardef\Deltag="7C01
\mathchardef\Phig="7C08
\mathchardef\Gammag="7C00
\mathchardef\Psig="7C09
\mathchardef\Lambdag="7C03
\mathchardef\Xig="7C04
\mathchardef\Pig="7C05
\mathchardef\Upsilong="7C07

\begin{document}

\title[Macdonald integrals and monodromy]{Macdonald integrals and monodromy}

\author{Jan Denef}
\address{University of Leuven, Department of Mathematics,
Celestijnenlaan 200B, 3001 Leu\-ven, Bel\-gium }
\email{ Jan.Denef@wis.kuleuven.ac.be}
\urladdr{http://www.wis.kuleuven.ac.be/wis/algebra/denef.html}

\author{Fran\c cois Loeser}

\address{École Normale Supérieure,
Département de math{\'e}matiques et applications,
45 rue d'Ulm,
75230 Paris Cedex 05, France
(UMR 8553 du CNRS)}
\email{Francois.Loeser@ens.fr}
\urladdr{http://www.dma.ens.fr/~loeser/}

%\date{Version of \today}
%\dedicatory{}

\begin{abstract}We prove several results on monodromies associated to
Macdonald integrals, that were used in our previous work on the
finite field analogue of a conjecture
of Macdonald. We also give a new proof 
of our formula expressing recursively
the zeta function of the local monodromy at
the origin of the discriminant of a finite Coxeter group in terms of
the degrees of the group. 
\end{abstract}

\maketitle

\section{Introduction}
\subsection{}Let $f :  \CC^n \rightarrow \CC$ be a polynomial map having a
singularity at the origin 0 in  $\CC^n$.
For $0 < \eta \ll \varepsilon \ll 1$,
the restriction of $f$ to 
$B (0, \varepsilon) \cap f^{-1} (D_{\eta} \setminus \{0\})$,
with $B (0, \varepsilon) $
the open ball of radius $\varepsilon$ centered at $0$ and
$D_{\eta}$ the open disk of radius $\eta$ centered at $0$,
is a locally trivial
fibration, the Milnor
fibration, onto 
$D_{\eta} \setminus \{0\}$ with fiber $F_{0}$, the Milnor fibre at $0$.
The action of a characteristic homeomorphism of this fibration on
cohomology
gives rise to 
the monodromy operator
$M : H^{\ast} (F_{0}, \QQ)
\rightarrow
H^{\ast} (F_{0}, \QQ)$.
Define 
the monodromy zeta function as
$$Z_{f, 0} (T): = \prod_{i \geq 0} [{\rm det} \, ({\rm Id} - TM, H^{i} (F_{0},
\QQ))]^{(-1)^{i + 1}}.
$$
In case the hypersurface $f = 0$ has an explicit embedded resolution,
a formula due to A'Campo \cite{Ac} may be used to compute
$Z_{f, 0}$. However, in general, calculating explicitely 
$Z_{f, 0}$ happens to be a quite dificult task.
\subsection{}Let $V$ be a complex vector space of finite dimension $n$
and
let $G$
be
a finite subgroup
of ${\rm GL} (V)$ generated by pseudo-reflections, {\it i.e.} endomorphisms
of finite order fixing pointwise an hyperplane. Such a group is called
a finite complex
reflection group.
Pseudo-reflections of
order 2 will be called reflections. Denote by $\CC [V]$
the algebra of polynomial functions on $V$.
By Chevalley's Theorem the ring of invariants $\CC [V]^G$ is a free
polynomial algebra on $n$ homogeneous invariant polynomials, whose
degrees $d_1$, \dots, $d_n$ only depend on $G$ and are called the
degrees of the group $G$.
For every pseudo-reflection in $G$ with corresponding hyperplane $H$, choose a
linear form $\ell_H$ defining $H$ and denote by $e (H)$
the order of the subgroup of elements of $G$ fixing $H$ pointwise.
Set $\Delta := \prod_H \ell_H^{e (H)}$, the product being over all 
pseudo-reflection hyperplanes in $G$. The induced function 
$\tilde \Delta : V / G \rightarrow \CC$  is the discriminant 
of $G$.
When
$V = \CC^n$, and
$G$ is  furthermore a subgroup of
${\rm GL} (\RR^n)$, $G$ is called a finite Coxeter group. In this case
the integers $e_H$ are all equal to 2.
Now consider $Z (T, G) :=
Z_{\tilde \Delta, 0} (T)$ the zeta function of the local monodromy of
the discriminant at the origin.

In the paper \cite{Indag} we proved the following remarkable recursion
formula
for $Z (T, G)$:

\begin{theorem}\label{deg}For $G$ a finite Coxeter group
we have
$$
\prod_{\cE \, \text{\rm connected  subgraph}}
Z (-T, G (\cE))^{(-1)^{\vert \cE \vert}}
=
\prod_{i = 1}^{n} \frac{1 - T^{d_{i}}}{1 -T}
$$
where the product on the left-hand side runs over all connected
subgraphs $\cE$ of the Coxeter diagram of $G$,
$G (\cE)$ denotes the Coxeter group with diagram $\cE$,
and $\vert \cE \vert$ the number of vertices of $\cE$.
\end{theorem}

The proof of Theorem \ref{deg} given in \cite{Indag} was based on some
new properties
of Springer's regular elements \cite{Springer}
in finite complex reflection groups, which have been since
further investigated by
Lehrer and Springer in \cite{LeSp1} and \cite{LeSp2}.
In fact, we computed in \cite{Indag} the zeta function of the local
monodromy of the discriminant for all irreducible
finite complex reflection groups. This involved a case by case
analysis, already for 
finite Coxeter groups.

\subsection{}\label{macdonot}Though
the techniques in the paper \cite{Indag} are mostly group theoretic, 
that work arouse from our study of a finite field analogue of 
Macdonald's conjecture.
Let us recall the statement of Macdonald's conjecture as 
formulated in \cite{Trans}.
Let $G$ be a finite subgroup of ${\rm GL} ({\RR}^n)$ which is 
generated by
reflections and let $q$
be a positive definite quadratic form
which is invariant under $G$.
We denote by $d_{1}$, \dots, $d_{n}$ the degrees of $G$.
Let $\ell_{1}$, \dots, $\ell_{N}$ be equations for the $N$ distinct 
reflection hyperplanes and set
$\Delta := (\prod_{i = 1}^{N} \ell_{i})^2$.
Macdonald's conjecture \cite{Macdo}, proved by Opdam \cite{Opdam},
is the following result.
\begin{theorem}The integral
\begin{equation}\label{defint}I (s)
:=\int_{\RR^n} \Delta (x)^s e^{-q(x)} dx
\end{equation}
may be expressed as
\begin{equation}\label{mdo}
I (s) =
\pi^{n/2} \kappa^s
\Bigl(
\prod_{i = 1}^n \frac{\Gamma (d_{i} s + 1)}{\Gamma (s + 1)}
\Bigr)
({\rm discr} \, q)^{- 1/2},
\end{equation}
with
\begin{equation}\label{kappa}
\kappa = \prod_{i = 1}^N \frac{q (\ell_{i})}{4},
\end{equation}
where we consider $\ell_{i}$ in $q (\ell_{i})$ as a vector in
$\RR^n$, identifying  $\RR^n$
with its dual, by means of the quadratic form $q$. 
\end{theorem}

\subsection{}Let us now state the finite field analogue.
Let $\FF $ be
a finite field of characteristic $p$
different from $2$. We consider 
a finite-dimensional $\FF$-vector space $V$,
a finite subgroup $G$  of ${\rm GL} (V)$
generated by reflections,
and $q$ a $G$-invariant, non degenerate,
symmetric bilinear form on $V$.
If $p$ does not divide $\vert G \vert$, one may define the degrees
of $G$, 
$d_{1}, \cdots, d_{n}$ as in the complex case. One also defines
$\Delta$ similarly.
Because
$p \not= 2$, we may define an element $\kappa$
of $\FF$ by (\ref{kappa}).
Fix a non trivial additive character
$\psi : \FF \rightarrow \CC$.
The analogue of the integral in (\ref{defint}) will be the character sum
\begin{equation}
S_{G} (\chi) := \sum_{x \in (U / G) (\FF)}
\chi (\Delta (x)) \psi (q (x))
\end{equation}
where
$\chi : \FF^{\times} \rightarrow \CC^{\times}$ is
a multiplicative character and $U$ denotes the complement of the
hypersurface
$\Delta (x) = 0$ in $V$. Here we write
$q (x)$ instead of $q (x, x)$.
The analogue of the Gamma function will be the Gauss sum $-g (\chi)$, where
$g (\chi) := -\sum_{x \in \FF^{\times}} \chi (x) \psi (x)$.
Our main result in the paper \cite{Trans} is the following finite
field analogue of
Macdonalds's conjecture.

\begin{theorem}\label{finite}Assume that
$p$ does not divide
$\vert G \vert$.  Then
$\kappa \not=0$ and,
for every multiplicative
character $\chi : \FF^{\times} \rightarrow \CC^{\times}$,
$$
S_{G} (\chi) =
(- 1)^{n} \phi ({\rm discr} q)
g (\phi)^{n} \phi (\kappa) \chi (\kappa) \prod_{i = 1}^{n}
\frac{g ((\phi \chi)^{d_{i}})}{g (\phi \chi)} \, ,
$$
where $\phi$ denotes the unique  multiplicative character of
order
2.
\end{theorem}

In the special case when $G$ is the symmetric group $S\sb n$, this
identity was conjectured by Evans \cite{E1}
and proved ten years later by him \cite{E2}
by using important work of 
Anderson \cite{Andsel}.

\subsection{}The aim of the present paper is doublefold. Firstly,
the proof of Theorem \ref{finite} in \cite{Trans}, not only used
Theorem \ref{deg}, but also some other results on the monodromy of
functions related to Macdonalds's integrals. These monodromy
calculations
were required in order to
use Laumon's product formula \cite{Lau}.
We present direct,
self-contained proofs of these results in
sections \ref{sec4}, \ref{sec5} and \ref{sec6}.
More precisely, in Theorem \ref{conn} of section \ref{sec4}
we express the global monodromy at the origin of the restriction of the
discriminant to the quotient of the quadric $q = 1$ by $G$ in terms of
local monodromies along strata of the discriminant.
In Theorem \ref{compl} of section 
\ref{sec5}, we prove a strange ``Complement formula'' expressing the 
local monodromy of the discriminant at the origin as
the sum of
the global monodromy at infinity of the restriction of the
discriminant to the quotient of the quadric $q = 1$ by $G$
and the global monodromy of the morphism
induced by $q^N$ at the origin in the quotient space
$U / G$.
In Theorem \ref{exactbar} of section \ref{sec6}, we explicitely compute 
the monodromy at the origin
of the morphism
induced by $q$  in the quotient space
$U / G$ and more generally for the constant sheaf replaced by the
pullback of a Kummer sheaf by the function defining the discriminant.
Secondly, we derive in section \ref{sec7} two consequences 
of Macdonald's formula (\ref{mdo}). The first one, Theorem \ref{max},
whose proof is
simple and elementary, is the calculation of the maximum of the function
$\Delta$ on the real points of the quadric defined by $q = 1$. This
formula is used in our proof of Theorem \ref{finite} and is the only
place in that proof relying on 
Macdonald's formula. Then, we explain how, using work of Anderson
\cite{And}
and Loeser and Sabbah \cite{LS} on determinants of Aomoto complexes
and determinants of integrals, one can derive Theorem \ref{deg} from
Macdonald's formula (in fact
Theorem \ref{deg} is equivalent to knowing the precise form of
the gamma factors in Macdonald's formula) making the full circle of the
story. 

At the end of the paper we give a complete list of the assertions in 
\cite{Trans}
whose proofs were postponed to the present work.

\subsection*{}{\small The present work was done around 1992, at about
the same time as the material in \cite{Indag} and \cite{Trans}. We
apologize for the long delay we have taken to finally write it down.
}

%\tableofcontents

\section{Coxeter arrangements}\label{sec2}

\subsection{Arrangements}Let
$V$ be a finite dimensional vector space over a
field
$k$ and $V \hookrightarrow \PP$
its canonical projective space compactification.
By an {\em hyperplane
arrangement} $\cA$ in $V$
we mean a finite set of affine
hyperplanes in $V$. Similarly, if $\PP$ is a projective space over $k$,
an hyperplane
arrangement in $\PP$ will be 
a finite set of projective
hyperplanes in $\PP$.

If $\cA$ is an arrangement in $V$,
we denote by 
$\overline \cA$ the projective arrangement in $\PP$
defined by taking the closure of the hyperplanes in $\cA$.
We shall denote by $V \setminus \cA$ the complement in
$V$ of the union of the hyperplanes 
belonging to $\cA$.
When $k = \RR$, we shall call a
connected component of
$V \setminus \cA$ a {\em chamber} of $V \setminus \cA$.
If all the hyperplanes contain $0$, the arrangement is said to
be {\em central}.
We call an endomorphism of $V$ a {\em reflection} if it has order
2 and fixes pointwise some hyperplane.

\begin{lem}\label{chamber}
Let $q$ be a non degenerate quadratic form on $\CC^n$ and denote by
$Q$ the quadric defined by $q = 1$ in $\CC^n$.
Let $\cA$ be a non empty
central arrangement in
$\RR^n$ such that the restriction
of $q$ to every non empty intersection of hyperplanes in $\cA$ 
is also non degenerate.Then 
the number of chambers ${\rm ch} (\RR^n \setminus \cA)$ of
$\RR^n \setminus \cA$ is equal to
$(- 1)^{n - 1} \chi (Q \setminus (\cA \otimes \CC))$.
Here $\cA \otimes \CC$ denotes the 
central arrangement in
$\CC^n$ obtained by extension of scalars.
\end{lem}

\begin{proof}When $|\cA| = 1$, the result is clear, since  
the Euler characteristic of the smooth quadric $Q$ in $\CC^n$ is equal to $1
+ (-1)^{n - 1}$. So assume 
$|\cA| > 1$ and
choose an hyperplane $H$ in $\cA$. Denote by
$\cA'$ the 
arrangement in
$\RR^n$ obtained by deleting $H$ from $\cA$ and
by $\cA''$
the 
arrangement in
$H$ obtained by intersecting $H$ with the hyperplanes in $\cA'$.
Since
${\rm ch} (\RR^n \setminus \cA) = 
{\rm ch} (\RR^n \setminus \cA')
+ {\rm ch} (H \setminus \cA'')$
and
$\chi (Q \setminus (\cA \otimes \CC)) = 
\chi (Q \setminus (\cA' \otimes \CC))
-
\chi ((Q \cap H) \setminus (\cA'' \otimes \CC))$,
the result follows by induction on the number of hyperplanes.
\end{proof}

\subsection{Coxeter arrangements}\label{ca}
We define a {\em classical Coxeter arrangement}
as a triple
$A = (V, G, q)$,
where $V$ is a finite dimensional vector space over $\RR$,
$G$ is a finite subgroup of ${\rm GL} (V)$
generated by reflections,
and $q$ is a $G$-invariant {\em positive definite}
symmetric bilinear form on $V$.
We define a {\em Coxeter arrangement over $\CC$}
as a triple
$A = (V, G, q)$,
where $V$ is a finite dimensional vector space over $\CC$,
$G$ is a finite subgroup of ${\rm GL} (V)$
generated by reflections,
and $q$ is a $G$-invariant 
symmetric bilinear form on $V$,
which arises by extension of scalars from a classical
Coxeter arrangement.
(In fact it follows from the argument given at the end of the proof
of Proposition 1.6 of 
\cite{Trans} that the last condition is always satisfied.)
We denote by $\cA_G$ the central arrangement consisting of all
reflection hyperplanes of $G$.

\begin{prop}\label{bc}Let $A = (V, G, q)$ be a
complex Coxeter arrangement
and denote by $Q$ the quadric defined by $q = 1$ in $V$.
Set $B = Q \cap (V\setminus \cA_G)$.
Then
$$
\chi (B) = (-1)^{n - 1} |G| \quad
\text{and} \quad \chi (B / G) = (-1)^{n - 1}.
$$
\end{prop}

\begin{proof}
Since a Coxeter group acts transitively on the corresponding set of
chambers (cf. \cite{Bbk} \S\kern .15em 3.1 Lemme 2), this follows 
from Lemma \ref{chamber}.
\end{proof}

\subsection{Canonical embedded resolution}Let
$V$ be a vector space of dimension $n$ over a
field
$k$ and let $\cA$ be an hyperplane
arrangement in $V$.
By an intersection space of $\cA$
we shall mean a non empty subset of
$V$ which is
the intersection of some non empty family of
hyperplanes in $\cA$.
By a stratum of $\cA$  we shall mean an intersection space of $\cA$ minus
the union of all hyperplanes of $\cA$ that do not contain that intersection 
space.
Similarly one defines intersection spaces
and strata of projective arrangements, and 
intersection spaces of 
$\overline \cA$ are just the closure of intersection spaces of $\cA$.
Let $\cA$ be an hyperplane arrangement in $V$ or 
in the canonical projective space compactification
$\PP$ of $V$ respectively.
We set $X_0 = V$ resp. $\PP$.
We define $h_1~: X_1 \rightarrow X_0$ to be the blowing up of all 
dimension zero intersection spaces of $\cA$ and then,
by induction, for
$2 \leq i \leq n - 1$, 
$h_i~: X_{i} \rightarrow X_{i - 1}$
to be the blowing up
of the union $Y_{i - 1}$ of the strict transforms of
all $i-1$-dimensional intersection spaces of $\cA$ 
in $X_{i - 1}$.
Note that $Y_{i - 1}$ is the disjoint union of these
strict transforms, hence is in particular smooth.
We set $X_{\cA} := X_{n - 1}$ and denote by $h_{\cA} : X \rightarrow X_0$
the composition of the morphisms $h_i$.
Remark that $X_{\cA}$ is smooth and that the inverse image of the union
of hyperplanes in $\cA$ by
$h_{\cA}$
is a divisor with (global) normal crossings, and that its  set of
irreducible
components 
is in natural bijection with the set of 
strata of $\cA$.

\subsection{}\label{resc}
Assume now $A = (V, G, q)$ is a 
complex Coxeter arrangement, fix
equations
$\ell_{1} = 0$, \dots, $\ell_{N}$= 0 for the $N$ distinct 
reflection hyperplanes and set
$\Delta := (\prod_{i = 1}^{N} \ell_{i})^2$.
Remark in this case $\cA_G$ and $\overline \cA_G$ have only one zero-dimension
stratum, namely the origin, hence $h_1$ is just the blow up of the origin.
It follows that
$h_{\overline \cA_G} : X_{\overline \cA_G} \rightarrow
\PP$ is
an embedded resolution of
the divisor $(\Delta = 0) \cup \overline
Q \cup H_{\infty}$ and the divisor $(\Delta = 0) \cup \overline
Q_0 \cup H_{\infty}$
with $\overline Q$, resp. $\overline Q_0$,
the closure of the locus of $q = 1$, resp. $q=
0$, in $\PP$ and 
$H_{\infty}$ the hyperplane at infinity. One should also remark that
the polar divisors of $\Delta$
and $q$ in 
$X_{\overline \cA_G}$ have only one irreducible component, namely
the strict transform of hyperplane at infinity $H_{\infty}$.
Furthermore, if we denote by $X_{{\overline \cA_G}, \overline Q}$ the strict transform
of
$Q$ in $X_{\overline \cA_G}$, the morphism
$X_{{\overline \cA_G}, \overline Q} \rightarrow  \overline Q$
yields an embedded resolution of the divisor 
$((\Delta = 0) \cup H_{\infty}) \cap \overline Q$
in $\overline Q$.

\section{Monodromy computations for group actions}\subsection{Monodromic Grothendieck group}\label{sec3}To agree with notation used in
\cite{Trans}, we shall use
the terminology
from
\cite{SGA 7} Exp. XIII - XIV concerning vanishing cycles. We shall explain
in \ref{zeta} how statements
about elements of the 
monodromic Grothendieck group
can be equivalently expressed in terms of the more classical
zeta functions of the monodromy used in
\cite{Milnor} \cite{Ac}.

We denote by $\bar \eta_{0}$ the generic geometric point of the
henselization of the complex affine line $\AA^{1}$ at $0$,
by $I_{0}$
its inertia group ({\it i.e.} the fundamental group of the complement of $0$
in a small disk around $0$) and by
$K_{I_{0}}$ the Grothendieck group
of
finite dimensional $\CC$-vector spaces with
$I_{0}$-action.
If $\cL$ is an object in $D^{b}_{c} (\GG_{m}, \CC)$,
the derived category of bounded complexes
of $\CC$-sheaves with constructible cohomology on the multiplicative group
$\GG_m := \AA^1 \setminus\{0\}$,
we denote by
$[\cL_{\bar \eta_{0}}]$ the class of
$\sum (- 1)^{i} [\cH^{i} (\cL)_{\bar \eta_{0}}]$
in $K_{I_{0}}$
and we set
$[\cL_{\bar \eta_{\infty}}] =
[{\rm inv}^{\ast}(\cL)_{\bar \eta_{0}}]$, where $\rm inv$
is the morphism $x \mapsto x^{-1}$.
If a finite group $G$ acts on $\cL$ we denote by
$[\cL^{G}_{\bar \eta_{0}}]$ the class of
$\sum (- 1)^{i} [\cH^{i} (\cL)^{G}_{\bar \eta_{0}}]$
and we define similarly $[\cL^{G}_{\bar \eta_{\infty}}]$.
For any character $\chi : I_{0} \rightarrow \CC^{\times}$ we denote
by $V_{\chi}$ the class in
$K_{I_{0}}$ of the rank one object with action given by $\chi$, hence
$V_{\chi} = [(\cL_{\chi})_{\bar \eta_{0}}]$,
with 
$\cL_{\chi}$ the local system of rank one on
$\GG_{m}$ and monodromy 
$\chi$ around the origin.
For any natural number $m \geq 1$,
we set
$V_{m} = [(\pi_{m \ast} \CC)_{\bar \eta_{0}}]$,
where $\pi_{m} : \GG_{m} \rightarrow \GG_{m}$ is given by $x \mapsto
x^{m}$,
so we have $V_{m} = \sum_{\chi^{m} = 1} V_{\chi}$.

\subsection{}\label{zeta}Let 
$\varrho$ be the standard topological generator
of $I_0$, corresponding to counter-clockwise rotation 
around $0$ in $\CC$. Since as an abelian group
$K_{I_{0}}$ is generated by the elements
$V_{\chi}$, there is a unique morphism of abelian groups
$$
Z : K_{I_{0}} \longrightarrow \CC (T)^{\times}
$$
sending $V_{\chi}$ to $Z (V_{\chi}) := (1 - T\chi (\varrho))^{- 1}$.
In particular we have
$Z (V_{m}) = (1 - T^m)^{- 1}$.
If $\cL$ is an object in $D^{b}_{c} (\GG_{m}, \CC)$,
$Z ([\cL_{\bar \eta_{0}}])$, resp.
$Z ([\cL_{\bar \eta_{\infty}}])$
is nothing else than the monodromy
zeta function
of $\cL$ at the origin, resp. at infinity,
of $\cL$. Remark that $Z$ is clearly an isomorphism onto its image.
Hence, depending on convenience, we shall either
formulate 
results in $K_{I_0}$ or in terms of monodromy
zeta functions.

\subsection{}\label{3.3}We consider the following geometric situation.
Let $j : U \rightarrow X$
be the immersion 
of a dense open subset $U$
in a smooth complex projective variety $X$.
We also assume to be given
a morphism
$f : U \rightarrow \GG_{m}$ and a proper morphism $g~: X
\rightarrow
\PP^{1}$
such that the diagram
\begin{equation*}\xymatrix{
U
\ar[d]_{f} \ar[r]^{j} & X \ar[d]^{g}\\
\GG_{m}
\ar[r]^{j_0}&\PP^{1}
}
\end{equation*}
is commutative, with $j_{0} : \GG_{m}\rightarrow \PP^{1}$ the
standard immersion.
We assume
$X-U$ is a divisor with normal crossings. We denote
by
$E_{i}, i\in J$, the irreducible components of
$X-U$. In particular the $E_{i}$'s are smooth. 
We set  $E=(g^{-1}(0))_{{\rm red}}= \cup_{i\in
J_{0}}E_{i}$, 
$E^{0}= E-\cup_{i\not\in J_{0}}E_{i}$, and  $g^{-1}(0)=
\sum_{i\in J_{0}}N_{i}E_{i}$.

\subsection{}We assume a finite group
$G$ acts on $X$ and that $U$ is stable.
We also assume the action of $G$ on $U$ is free and that $f$ is
$G$-invariant.
For $s\in E$ we denote by $G_{s}$ the stabilizer at $s$ and
by $[(R \psi_{g}\CC)_{s}^{G_{s}}]$ the class of
$\Sigma(-1)^{j}[(R^{j}\psi_{g}\CC)^{G_{s}}_{s}]$ in
$K_{I_{0}}$, with $R \psi_{g}\CC$
the complex of nearby cycles of $g$.
For $s\in E^{0}$, we set $J^{s}=\{ i\in J~; s\in
E_{i}\}\subset J_{0}$ and we denote by
$C_s$ the set of connected components of the Milnor fiber of $g$ at $s$.
There is a natural action of $G_s$ on $C_s$.

The following result is proved in 
\S\kern .15em 3
of \cite{BSMF}:

\begin{prop}\label{mono}
Assume the previous assumptions hold.
\begin{enumerate}
\item[(1)]Let $\cS$ be a partition of $E^{0}$ into
constructible subsets $S$
such that $|G_{s}|$ and $[(R \psi_{g}\CC)^{G_{s}}_{s}]$ are
constant on  every $S$. Then the following holds in
$K_{I_{0}}$:
\begin{equation*}
[( R f_{!} \CC)^{G}_{\bar\eta_{0}}]= |G|^{-1}\sum_{S\in\cS}\chi
(S,\CC)|G_{s}| [(R\psi_{g}\CC)^{G_{s}}_{s}],
\end{equation*}
where $s \in \CC$. 
\item[(2)]For every closed point $s$ of $E$,
the stabilizer $G_s$ acts freely on  $C_s$.
\item[(3)]Assume for every $i\in J_{0}$ and
every $\sigma \in G$, $\sigma (E_{i})=E_{i}$ or $\sigma (E_{i})\cap
E _{i}=\emptyset$.
Then, for every $s\in E^{0}$,
\begin{equation*}
 [(R \psi_{g}\CC)^{G_{s}}_{s}] 
 = 
\begin{cases}0 &
\text{if $ |J^{s}| > 1$},\\
(\CC^{C_s/G_s})^{\vee} & \text{if $|J^{s}| = 1$},
\end{cases}
\end{equation*}
with $(\CC^{C_s/G_s})^{\vee}$ the dual of the module
$\CC^{C_s/G_s}$ endowed with its natural $I_{0}$-action.
\end{enumerate}
\end{prop}

\subsection{}\label{nd}In fact, we shall also need to consider the case where
the morphism $f$,
instead of having an extension to an actual
morphism $g~: X
\rightarrow
\PP^{1}$,
only extends to a rational map
$g~: X
\rightarrow
\PP^{1}$.
Denote by $g^{-1} (0)$ and $g^{-1} (\infty)$
the divisor of zeroes and poles respectively.
Remark that $g$ is a morphism outside
$g^{-1} (0) \cap g^{-1} (\infty)$ and write
$g^{-1} (\infty) = \sum_{i \in J_{\infty}} N_i E_i$.
We replace
the condition ``$X \setminus U$ is a divisor with normal
crossings''
by
``$X \setminus U$ is a divisor with normal
crossings on a Zariski neighbourhood of $g^{-1} (0)$''.

The following generalisation of Proposition \ref{mono}
is also proved in 
\S\kern .15em 3
of \cite{BSMF} and shows that the exceptional divisors of the
additional blowing ups needed to make the extensions of $f$ and $g$
regular
do not contribute to our monodromy calculations:

\begin{prop}\label{monog}Assume for every $i$ in 
$J_0 \cup J_{\infty}$ and every $\sigma$ in $G$, 
$\sigma (E_i) = E_i$ or
$\sigma (E_i) \cap E_i = \emptyset$. Then the statements in 
Proposition \ref{mono} still hold in the  setting \ref{nd}.
\end{prop}

\section{Expressing the global monodromy at the origin in term of
local monodromies}\label{sec4}

\subsection{}\label{coxnot}Let $A = (V, G, q)$ be a complex 
Coxeter arrangement.
As in \ref{ca} we fix equations for the hyperplanes and consider the
corresponding function $\Delta$.
Let
$(R\psi_{\Delta} (\CC))_{0}$ be the stalk  at zero of the complex
of nearby cycles with respect to $\Delta$. We
set
$M_{G} := [(R\psi_{\Delta} (\CC))_{0}^G]$ and
$\bar M_{G} := (-1)^{n - 1} [(R\psi_{\Delta} (\CC))_{0}^G]$,
with $n$ the dimension of $V$.
We also set $B := Q \setminus \cA_G$.
\begin{theorem}\label{conn}The relation 
$$[(R \Delta_{\vert B !} \CC)^{G}_{\bar \eta_{0}}]
=
(- 1)^n \sum_{{\cE \,
\text{\rm connected \, subgraph}}\atop{G (\cE) \not= G}} \bar M_{G
(\cE)}$$
holds in
$K_{I_{0}}$.
\end{theorem}

\begin{proof}Note that
$M_G = [(R \Delta_{!} \CC)^{G}_{\bar \eta_{0}}]$, since $\Delta$ is
homogeneous. Thus we can calculate
$M_G$ by applying Proposition \ref{monog}
to the resolution
$X_{{\overline \cA_G}} \rightarrow  \PP$ defined
in \ref{resc}. 
Since the Euler characteristic of a complex algebraic variety
with a free $\GG_m$-action is zero (see, {\it e.g.}, \cite{bb}),
we have to sum in the formula of \ref{mono} (1) for $M_G$ only over
strata lying inside the strict transform in 
$X_{{\overline \cA_G}}$ of the exceptional divisor $H$
of the blow up $h_1 : X_1 \rightarrow V$ of $0$ in $V$.
Let $Z$ be the complement in $H$ of the
strict transform
in $X_1$ of the locus of $\Delta = 0$, and let $\cal W$ 
be a partition of $Z$ into constructible subsets $W$
on which $G_s$ and
$[(R\psi_{\Delta \circ h_1} (\CC))_{s}^{G_s}]$
are constant. Then the above discussion yields
$$
M_G = 
\vert G \vert^{-1}
\sum_{W \in \cW}
\chi (W) \vert G_s \vert [(R\psi_{\Delta \circ h_1} (\CC))_{s}^{G_s}],
$$
where $s$ is any point in $W$. A similar formula holds also
for 
$M_{G_S}$, where $G_S$ is defined below.
Applying now Proposition \ref{monog} to the resolution
$X_{{\overline \cA_G}, \overline Q} \rightarrow  \overline Q$ defined
in \ref{resc}, and using \cite{Bbk} \S\kern .15em 3.3 Proposition 1,
we obtain
the relation
$$
[(R \Delta_{\vert B !} \CC)^{G}_{\bar \eta_{0}}]
=
\frac{1}{|G|}
\sum_{S \text{\rm stratum}}
\chi (S \cap Q) \, |G_S| \, M_{G_S},
$$
where, for any subset $S$ of $V$, $G_S$ is the group generated by 
the reflections with respect to the walls containing $S$.
Here by stratum, we mean non dense stratum of the arrangement
associated to $G$.
Using Lemma \ref{chamber}, we obtain
\begin{equation*}
\begin{split}
[(R \Delta_{\vert B !} \CC)^{G}_{\bar \eta_{0}}]
&=
\frac{1}{|G|}
\sum_{S \text{\rm stratum}, S \not= \{0\}}
(-1)^{{\rm dim} S  - 1}
(\sum_{{F \, \text{\rm face}} \atop {F \subset S}} 1)
 |G_S| M_{G_S}\\
&=
\frac{1}{|G|}
\sum_{F \, \text{\rm face}, F \not= \{0\}}(-1)^{{\rm dim} F  - 1}
|G_F| M_{G_F},\\
\end{split}
\end{equation*}
where by ``face'' we mean any
face of codimension $\geq 1$ of some chamber
of the arrangement associated to $G$. The faces which are contained in
a stratum $S$ are indeed exactly the chambers of the set of real points of $S$.
Let $C_0$ be a fixed chamber. For each face $F$, there exists a
unique
face $F_0$ of $C_0$ such that there exists $w$ in $G$ with
$F_0 = w (F)$ (see \cite{Bbk} \S\kern .15em 3.3, Remarque 1).
We say that the faces $F$ and $F_0$ are {\em related}.
We have
\begin{equation*}
\begin{split}
[(R \Delta_{\vert B !} \CC)^{G}_{\bar \eta_{0}}]
&=
\frac{1}{|G|}
\sum_{{F_0 \, \text{\rm face of} \, C_0}\atop {F_0 \not= \{0\}}}
\sum_{{F \, \text{\rm face}} \atop  {F \, \text{\rm related to} \, F_0}}
(-1)^{{\rm dim} F  - 1}
|G_F| M_{G_F}\\
&=
\frac{1}{|G|}
\sum_{{F_0 \, \text{\rm face of} \, C_0}\atop {F_0 \not= \{0\}}}
(-1)^{{\rm dim} F_0  - 1}
|G_{F_0}| M_{G_{F_0}} N (F_0),\\
\end{split}
\end{equation*}
with
$N (F_0)$ the number of faces related to $F_0$.
Note that 
$N (F_0)$ is equal to the total number of chambers divided by the
number of chambers in the arrangement associated to $G_{F_0}$,
since every chamber has exactly one face related to $F_0$, and the
number of chambers containing the same face $F$ is equal 
to
the number of 
chambers in the arrangement associated to $G_{F}$.
Hence, since the number of chambers of a Coxeter arrangement is
equal to the order of the group, we obtain
$$
[(R \Delta_{\vert B !} \CC)^{G}_{\bar \eta_{0}}]
=
\sum_{{F_0 \, \text{\rm face of} \, C_0}\atop {F_0 \not= \{0\}}}
(-1)^{{\rm dim} F_0  - 1}
M_{G_{F_0}}.
$$
Associating to $F_0$ the subgraph of the Coxeter diagram $\cG$
of $G$
whose vertices correspond to the walls of $C_0$ that contain
$F_0$, we get
$$[(R \Delta_{\vert B !} \CC)^{G}_{\bar \eta_{0}}]
=
(- 1)^n \sum_{{\cE \,
\text{\rm proper subgraph of}}\, \cG}  \bar M_{G (\cE)}.$$ 
Since, by Corollary 3.3 of \cite{Indag},
which is given in
\cite{Indag}
a direct
self contained
proof,
$\bar M_{G
(\cE)} = 0$ when $\cE$ is not connected, the result follows.
\end{proof}

\section{The formula of the complement}\label{sec5}
\begin{lem}\label{lem1}The following holds in
$K_{I_{0}}$:
$$
M_G = [(R \Delta_{\vert B_0 !} \CC)^{G}_{\bar \eta_{0}}] + 
[(Rq^N_{|U !} \CC)^G_{\bar \eta_0}].
$$
\end{lem}

\begin{proof}For any element $W$ of $K_{I_0}$,
the monodromy zeta function $Z (W)$ of $W$
is equal to $\exp \sum_{j = 1}^{\infty}
{\rm Tr} \,(\varrho^j, W) T^j/j$,
where $\varrho$ is as in \ref{zeta} (see, {\it e.g.}
\cite{Milnor}, p.77).
Hence, we have to prove that
\begin{equation}\label{tr}
{\rm Tr} \, (\varrho^j, M_G) = 
{\rm Tr}\,(\varrho^j, (R \Delta_{\vert B_0 !} \CC)^{G}_{\bar \eta_{0}}) +
{\rm Tr}\,(\varrho^j, (Rq^N_{|U !} \CC)^G_{\bar \eta_0}),
\end{equation}
for every $j$ in $\NN \setminus \{0\}$.
Since both functions $\Delta$ and $g^N$ are homogeneous
of degree $2N$, the map $h : \CC^n \rightarrow \CC^n$, $x \mapsto
e^{2 i \pi/2N} x$ induces the monodromy action $\varrho$
on the elements of 
$K_{I_0}$ appearing in (\ref{tr}).
Hence the traces in (\ref{tr}) are equal to
the Euler characteristic of the fixed point manifold
of $h^j$ restricted to the generic fibre of respectively
$\Delta_{\vert U/G}$,
$\Delta_{\vert B_0/G}$
and $q^N_{\vert U/G}$ (see, {\it e.g.}, Lemma 9.5
of \cite{Milnor}).
Thus,
we may write
\begin{equation*}
\begin{split}
{\rm Tr} \, (\varrho^j, M_G) =&
\chi \Bigl(\bigl\{x \in U/G \bigm \vert
\Delta (x) = 1, h^j (x) = x \mod G
\bigr\}
\Bigr)
\\=&
\chi \Bigl(\bigl\{x \in U/G \bigm \vert
q (x) = 0, \Delta (x) = 1, h^j (x) = x \mod G
\bigr\}\Bigr)\\
&+
\chi \Bigl(\bigl\{x \in U/G \bigm \vert
q (x) \not=0, \Delta (x) = 1, h^j (x) = x \mod G
\bigr\}\Bigr)\\
\end{split}
\end{equation*}
and
\begin{equation*}
\begin{split}
{\rm Tr}\,(\varrho^j, (R \Delta_{\vert B_0 !} \CC)^{G}_{\bar \eta_{0}}) + &
{\rm Tr}\,(\varrho^j, (Rq^N_{|U !} \CC)^G_{\bar \eta_0})
=
{\rm Tr}\,(\varrho^j, (R \Delta_{\vert B_0 !} \CC)^{G}_{\bar \eta_{0}})\\
&+
\chi \Bigl(\bigl\{x \in U/G \bigm \vert
q^N (x) =1, \Delta (x) \not= 0, h^j (x) = x \mod G
\bigr\}\Bigr).\\
\end{split}
\end{equation*}
Hence it is enough to prove that
\begin{equation}\label{5.1.2}
\begin{split}
\chi \Bigl(\bigl\{x \in U/G \bigm \vert
q (x) \not=0, &\, \Delta (x) = 1, h^j (x) = x \mod G
\bigr\}\Bigr)
=\\
\chi \Bigl(\bigl\{x \in U/G & \bigm \vert
q^N (x) =1, \Delta (x) \not= 0, h^j (x) = x \mod G
\bigr\}\Bigr).
\end{split}
\end{equation}
This follows from the fact that the variety
$$
\bigl\{(x, z) \in U/G \times \CC^{\times} \bigm \vert
q (x) \not=0, \Delta (x) = 1, h^j (x) = x \mod G, z^{2N} = q^N (x)
\bigr\}
$$
is isomorphic to
the variety
$$
\bigl\{(x, z) \in U/G \times \CC^{\times} \bigm \vert
q^N (x) =1, \Delta (x) \not=0, h^j (x) = x \mod G, z^{2N} = \Delta (x)
\bigr\}
$$
through the morphism $(x, z) \mapsto (xz^{-1}, z^{-1})$,
whose inverse is given by the same formula, since the two varieties are
étale covers of degree $2N$ of the two varieties occuring in
(\ref{5.1.2}).\end{proof}

\begin{prop}\label{three}
We have
$$
[(R \Delta_{\vert B_0 !} \CC)^{G}_{\bar \eta_{0}}]
=
[(R \Delta_{\vert B_0 !} \CC)^{G}_{\bar \eta_{\infty}}]
=
[(R \Delta_{\vert B !} \CC)^{G}_{\bar \eta_{\infty}}]
$$
in $K_{I_0}$.
\end{prop}

\begin{proof}
Here again 
we apply Proposition \ref{monog} to the resolution
$X_{{\overline \cA_G}, \overline Q} \rightarrow  \overline Q$, defined
in \ref{resc}, to compute 
$[(R \Delta_{\vert B !} \CC)^{G}_{\bar \eta_{\infty}}]$.
And similarly, replacing $Q$ by $Q_0$, one computes
$[(R \Delta_{\vert B_0 !} \CC)^{G}_{\bar \eta_{\infty}}]$.
Direct observation shows that both are equal.
To prove the equality
$[(R \Delta_{\vert B_0 !} \CC)^{G}_{\bar \eta_{0}}]
=
[(R \Delta_{\vert B_0 !} \CC)^{G}_{\bar \eta_{\infty}}]$,
one remarks  that 
$\Delta$ induces a $G$-equivariant fibration
$B_0 (\CC) \rightarrow \CC^{\times}$, hence
$(R \Delta_{\vert B_0 !} \CC)_{\bar \eta_{0}}$
and
$(R \Delta_{\vert B_0 !} \CC)_{\bar \eta_{\infty}}$
are already isomorphic as complexes of sheaves with $G$-action.
\end{proof}

\begin{remark}It seems quite likely
that
$(R \Delta_{\vert B_0 !} \CC)_{\bar \eta_{\infty}}$
and
$(R \Delta_{\vert B !} \CC)_{\bar \eta_{\infty}}$
are already isomorphic as complexes of sheaves with $G$-action.
\end{remark}

Now we are able to deduce the following result.
\begin{theorem}[Formula of the complement]\label{compl}
The equality
$$
[(R \Delta_{\vert B !} \CC)^{G}_{\bar \eta_{\infty}}]
+
[(Rq^N_{|U !} \CC)^G_{\bar \eta_0}]
=
M_G
$$
holds in $K_{I_0}$.
\end{theorem}

\begin{proof}Follows directly from Lemma
\ref{lem1} and Proposition \ref{three}.
\end{proof}

\section{Calculation of
$[(Rq_{|U !} \Delta^{\ast}\cL_{\chi})^G_{\bar \eta_0}]$}\label{sec6}
We first begin by proving
the following general result.
\begin{prop}\label{top}Let $f$ in
${\RR} [x_1, \ldots, x_n]$ be a homogeneous polynomial of degree
$N$.
Let $G$ be a finite subgroup of ${\rm GL} ({\RR}^n)$, and let $q$
be a positive definite quadratic form on 
${\RR}^n$ which is $G$-invariant.
We assume that $f \circ \sigma = {\rm det} (\sigma) f$
for every $\sigma$ in $G$.
We set $\Delta = f^2$, we  denote by  $U$
the complement in  ${\CC}^n$ of the hypersurface
$\Delta = 0$, and by $Q$ the quadric $q = 1$ in ${\CC}^n$. We 
set $B := Q \cap U$ and we denote by 
$\phi$ the unique character of order 2 of $I_0$.
Assume that 
\begin{equation}\label{star}
H^i_c (B, \Delta^{\ast}\cL_{\chi}) = 0
\quad \hbox{\rm for} \quad i \not= n- 1
\end{equation}
and
\begin{equation}\label{2star}
{\rm dim} \, H^{n - 1}_c (B, \Delta^{\ast}\cL_{\chi})^G = 1,
\end{equation}
for almost
all characters $\chi$ of finite order of $I_0$.
Then, for almost all such 
$\chi$, we have
\begin{equation}\label{restar}
[(Rq_{|U !} \Delta^{\ast}\cL_{\chi})^G_{\bar \eta_0}]
= 
(-1)^{n - 1} V_{\phi^{n + N} \chi^N}
\end{equation}
in the Grothendieck group $K_{I_{0}}$.
Here ``almost all'' means ``outside a finite set''.
\end{prop}

\begin{proof}Denote by 
$(Rq_{|U !}\Delta^{\ast}\cL_{\phi \chi})^{\rm det}$
the part of 
$(Rq_{|U !}\Delta^{\ast}\cL_{\phi \chi})$
on which the $G$-action is given by multiplication by the determinant.
By a direct  adaptation of the proof of Lemma 2.3.1 (1) of
\cite{Trans} we get an isomorphism
\begin{equation}\label{6.1.4}
(Rq_{|U !} \Delta^{\ast}\cL_{\chi})^G \simeq 
(Rq_{|U !}\Delta^{\ast}\cL_{\phi \chi})^{\rm det}.
\end{equation}
Thus, it suffices to show that
$$[(Rq_{|U !} \Delta^{\ast}\cL_{\chi})^{\rm det}_{\bar \eta_0}]
=
(-1)^{n - 1} V_{\phi^{n} \chi^N}.
$$
Now consider the compact set $Q ({\RR}) = Q \cap {\RR}^n$
and set $C := Q ({\RR})\cap U$.
Since $Q ({\RR})$ is compact and 
$\cL_{\chi}$ is constant on 
$\Delta (C)$,
the set $C$ determines a cycle class
$[C]_{\chi}$ in
$$
H^{\Delta-{\rm proper}}_{n - 1} (B, \Delta^{\ast}\cL_{\chi})
\simeq
H_{\Delta-{\rm proper}}^{n - 1} (B, \Delta^{\ast}\cL_{\chi})
$$
which is non zero since
$$
\int_C \Delta^s \frac{dx}{dq} \not= 0
$$
for every $s >0$.
We remark that, for $\sigma$ in $G$, we have
$$
\sigma ([C]_{\chi}) = {\rm det} (\sigma) [C]_{\chi},
$$
because $[C]_{\chi}$ is induced by an element of
$H_{n - 1} (Q (\RR), \RR)$ on which $G$ acts 
by multiplication with the  determinant.
For almost all
$\chi$,
the canonical morphism
$$
H^{n - 1}_c (B, \Delta^{\ast}\cL_{\chi})^{\rm det}
\longrightarrow
H^{n - 1}_{\Delta-{\rm proper}} (B, \Delta^{\ast}\cL_{\chi})^{\rm det}
$$
is an isomorphism by  Proposition 4.2.7 of \cite{LS}.
Note that $H^{n - 1}_c (B, \Delta^{\ast}\cL_{\chi})^{\rm det}$
has rank one, for almost all $\chi$, because of (\ref{2star})
and (\ref{6.1.4}).
Thus, when $\chi$ is general enough, the cycle class $[C]_{\chi}$
is a generator of $H^{n - 1}_c (B, \Delta^{\ast}\cL_{\chi})^{\rm det}$.
Choose a topological generator  $\varrho$ of $I_0$.
It is enough to prove that
$$
\varrho ([C]_{\chi}) = (-1)^n \chi^N (\varrho)[C]_{\chi}.
$$
Now remark that the map $x \mapsto \exp({2 \pi i \theta/2}) x$,
with $\theta \in [0,1]$ is a realization of the monodromy of $q$ which
is induced by $- {\rm Id}_{\CC^n}$.
Let  $\chi$ be the character sending
$\varrho$ to $\exp ({2 \pi i a/k})$.
Since $\Delta (\exp ({2 \pi i \theta /2}) x) = 
\exp({2 \pi i N \theta}) (\Delta (x))$, we obtain that
\begin{equation*}
\begin{split}
\varrho ([C]_{\chi})  &= \exp ({2 \pi i a N/k}) {\rm det}
(-{\rm Id}_{\CC^n}) [C]_{\chi}\\
&= \chi^N (\varrho)(-1)^n [C]_{\chi},\\
\end{split}
\end{equation*}
and the result follows.
\end{proof}

\begin{prop}\label{ab}Assume we are in the Coxeter 
setting \ref{coxnot}.
There exists integers $\bar a$ and 
$\bar b$ satisfying $\bar a + \bar b = (-1)^{n - 1}$ ,
such that the following relations hold:
\begin{enumerate}
\item[(1)]$[(Rq_{|U !} \CC)^G_{\bar \eta_0}] = \bar a V_1 + \bar b
  V_{\phi} = (\bar a - \bar b)V_1
  +
\bar b
  V_{2}$.
\item[(2)]For every character $\chi$ of $I_0$, we have
$[(Rq_{|U !} \Delta^{\ast}\cL_{\chi})^G_{\bar \eta_0}]
=
\bar a V_{\chi^N} + \bar b V_{\phi \chi^N}.
$
\item[(3)]$[(Rq^N_{|U !} \CC)^G_{\bar \eta_0}] = (\bar a - \bar b)V_N
  +
\bar b
  V_{2N}$.
\end{enumerate}
\end{prop}

\begin{proof}Since $q$ is homogeneous of degree 2 and $\Delta$ is homogeneous,
the map $- {\rm Id}_{\CC^n}$ induces the monodromy action
on 
$(Rq_{|U !} \CC)^G_{\bar \eta_0}$.
This proves the existence of
$\bar a$ and $\bar b$ in (1).
That 
$\bar a + \bar b = (-1)^{n - 1}$
follows directly form Proposition \ref{bc}, since the virtual rank of
$[(Rq_{|U !} \CC)^G_{\bar \eta_0}]$ equals $\chi (B) /\vert G\vert$.
To prove (2) we first show that
\begin{equation}\label{jjj}
[(Rq_{|U !} \Delta^{\ast}\cL_{\chi})^G_{\bar \eta_0}]
=
[(Rq_{|U !} q^{N \ast}\cL_{\chi})^G_{\bar \eta_0}].
\end{equation}
To show (\ref{jjj}), we use the fact that
\ref{mono} (1), and 
\ref{mono} (3) for $\vert J_s \vert > 1$, remain valid when the 
constant sheaf  is replaced by
a local system on $U$ with $G$-action (see Remarque 3.4.2 of
\cite{BSMF}). Moreover, this remains valid in the more general situation
of Proposition \ref{monog}.
We apply this to the rational map
$g : X_{\overline \cA_G} \rightarrow \PP^1$
induced by $q$.
In
this case, using the notation of \ref{3.3},
we have $E = \widetilde H \cup \widetilde Q$
where
$\widetilde H$ is the strict transform
in $X_{\overline \cA_G}$ of the exceptional
divisor $H$ of the blow up 
$h_1 : X_1 \rightarrow V$ of $0$ in $V$,
and $\widetilde Q$ is the 
strict transform
in $X_{\overline \cA_G}$ of the locus of $q = 0$ in $V$.
Since the Euler characteristic of a complex algebraic variety
with a free $\GG_m$-action is zero (see, {\it e.g.}, \cite{bb}),
we may sum in \ref{mono} (1) only over strata lying inside
$\widetilde H$.
In fact, because of \ref{mono} (3), we may even only sum
over strata lying inside
$\widetilde H \setminus \bigcup_{E_i \not= \widetilde H}E_i
\simeq H \setminus S$,
where $S$ is the
strict transform in $X_1$ of the locus of $\Delta q = 0$.
Equality (\ref{jjj}) follows now directly from the above discussion
and the fact that $h_1^{\ast} q^{N \ast} \cL_{\chi}$
and
$h_1^{\ast} \Delta^{\ast} \cL_{\chi}$
are locally isomorphic, as local systems with $G$-action,
on a neighbourhood
of $H \setminus S$ in $X_1$.

But $q^{N \ast}\cL_{\chi} \simeq q^{\ast}\cL_{\chi^{N}}$,
hence
$$[(Rq_{|U !} q^{N \ast}\cL_{\chi})^G_{\bar \eta_0}]
=
[(Rq_{|U !} (\CC) \otimes \cL_{\chi})^G_{\bar \eta_0}]
$$
by the projection formula,
and
$$
[(Rq_{|U !} \Delta^{\ast}\cL_{\chi})^G_{\bar \eta_0}]
= 
[(Rq_{|U !} \CC)^G_{\bar \eta_0}] \otimes V_{\chi},$$
which shows that (2) follows from (1).
Since $q^N = \pi_N \circ q$, (3) follows from (1).
\end{proof}

We now determine the exact value of $\bar a$ and $\bar b$.
\begin{theorem}\label{exactbar}Assume we are in the Coxeter 
setting \ref{coxnot}. Then, for every
character 
$\chi$ of $I_0$, we have
\begin{equation*}\label{3star}
[(Rq_{|U !} \Delta^{\ast}\cL_{\chi})^G_{\bar \eta_0}]
= 
(-1)^{n - 1} V_{\phi^{n + N} \chi^N}
\end{equation*}
in the Grothendieck group $K_{I_{0}}$.
\end{theorem}

\begin{proof}Let us check that 
conditions (\ref{star}) and (\ref{2star}) are verified for almost all
$\chi$.
Consider the open immersion
$j : B \hookrightarrow X_{{\overline \cA_G}, \overline Q}$.
The canonical morphism
$$
Rj_! ((\Delta^{\ast} \cL_{\chi})_{|B})
\longrightarrow 
Rj_{\ast} ((\Delta^{\ast} \cL_{\chi})_{|B})
$$
is an isomorphism for almost all $\chi$,
hence the canonical morphism
$$H^i_c (B, \Delta^{\ast}\cL_{\chi})
\longrightarrow 
H^i (B, \Delta^{\ast}\cL_{\chi})
$$
is an isomorphism for almost all $\chi$.
Since $B$ is affine, $H^i (B, \Delta^{\ast}\cL_{\chi})$ is zero for $i
> n-1$,
hence, by Poincar{\'e} duality, it follows that (\ref{star})
is verified for almost all
$\chi$. For such a $\chi$, the rank of
$H^{n - 1}_c (B, \Delta^{\ast}\cL_{\chi})^G$ is equal $(-1)^{n - 1}$
times  the Euler
characteristic of $B / G$, so (\ref{2star}) follows from
Proposition \ref{bc}.
Now the result follows
by putting together Proposition \ref{top} and
Proposition \ref{ab}.
\end{proof}

\section{From Macdonald integrals to monodromy}\label{sec7}
\subsection{}In this section we shall
work in the framework of \ref{macdonot}.

Set $S := \{x \in \RR^n \, | \, q (x) = 1\}$ and 
observe 
that
\begin{equation}
\sqrt{q (\ell_{i})} = {\rm Max}_{x \in S}
\ell_{i}(x) \, . 
\end{equation}

Before explaining the relations with monodromy,
let us first derive the following 
interesting consequence of Macdonald's formula (\ref{mdo}),
which is used in the paper \cite{Trans}. We do not know a direct proof
of this result.

\begin{theorem}\label{max}
We have
$$
{\rm Max}_{x \in S} \Delta (x) =
\kappa \, \frac{\prod_{i = 1}^n d_{i}^{d_{i}}}{N^N}.
$$
\end{theorem}

\begin{proof}We have
\begin{equation}
{\rm Max}_{x \in S} \Delta (x) =
\lim_{s \rightarrow + \infty} \Bigl(
\int_{S} \Delta (x)^s \frac{\vert dx\vert}{\vert dq\vert}
\Bigr)^{1/s}
\end{equation}
and
\begin{equation}
\int_{S} \Delta (x)^s \frac{\vert dx\vert}{\vert dq\vert}
= 
\frac{1}{\Gamma (Ns + \frac{n}{2})}
\int_{\RR^n} \Delta (x)^s e^{-q(x)} dx.
\end{equation}
Hence we deduce from (\ref{mdo}), that
\begin{equation}
{\rm Max}_{x \in S} \Delta (x) = \kappa \,
\lim_{s \rightarrow + \infty} \Bigl(
\frac{1}{\Gamma (Ns + \frac{n}{2})}
\prod_{i = 1}^n \frac{\Gamma (d_{i} s + 1)}{\Gamma (s + 1)}
\Bigr)^{1/s}.
\end{equation}
The result follows now from
Stirling's formula
$\Gamma (x + 1) \simeq  \sqrt{2 \pi} x^{x + 1/2} e^{- x}$
and the relation $\sum_{i = 1}^{n} (d_{i} - 1) = N$.
\end{proof}

\begin{remark}\label{last}Theorem 3.3 in \cite{Trans} follows
directly from Theorem \ref{max},
since the morphism
$\Delta_{\vert B} : B \rightarrow \GG_m$
has a compactification whose restriction
to the locus at infinity
is analytically trivial locally at each point of that locus, and
because
$\Delta_{\vert B}$ has only non degenerate
critical points, that are all  conjugate under $G$.
As compactification one takes for instance
the projection onto $\GG_m$
from the closure in
$X_{\overline \cA_G, \overline Q} \times \GG_m$ of the graph of
$\Delta_{\vert B}$ (see section \ref{resc}).
For this it is essential that the support of the divisor 
of $\Delta$ on 
$X_{\overline \cA_G, \overline Q}$ is exactly 
$X_{\overline \cA_G, \overline Q} \setminus B$.
To prove the assertion about the 
critical points of 
$\Delta_{\vert B}$, it is enough to show that
$\Delta_{\vert B}$ has only isolated critical points, because
then $(-1)^{n- 1} \chi (B)$ is equal to the sum of the Milnor numbers
of the critical points of $\Delta_{\vert B}$.
But $|G| = (-1)^{n - 1} \chi (B)$
by Proposition \ref{bc}, and there are at least
$|G|$ critical points because the action of $G$ on $B$ is free.
To see that $\Delta_{\vert B}$ has only isolated critical points, note that
the zero locus $Z$ of the section $d \Delta /\Delta$
of $\Omega^1_{X_{\overline \cA_G, \overline Q}}
({\rm log} (X_{\overline \cA_G, \overline Q} \setminus B))$
is closed in the proper variety $X_{\overline \cA_G, \overline Q}$,
but contained in the affine variety $B$, hence $Z$ is just a finite set.
\end{remark}

\subsection{}In their paper \cite{LS}, Loeser and Sabbah
gave a general formula for computing the determinant of a matrix
whose entries are integrals of algebraic differential forms multiplied
by a product of complex powers of polynomials. The
gamma factors that appear in this formula
are described in terms of monodromies associated with
the family of polynomials. The proof involved the computation of the
determinant of a complex of twisted 
differential forms introduced by K. Aomoto. The computation of such a
determinant of a complex of twisted 
differential forms was carried out independently (and slightly before) by
Anderson in
\cite{And}.
Since stating the general formula would lead us to far away from the
core of the present work, we shall only quote the result in 
\cite{LS} in a very special case, which will be enough for what we need in
this paper.

\begin{theorem}\label{lsa}Let $X$ be a smooth connected
algebraic variety over $\RR$ of dimension $n$, and let $f : X
\rightarrow \GG_{m, \RR}$ be a morphism of real 
algebraic varieties.
Let $\omega$ be a global section of $\Omega_X^n$ and let $\Gamma$ be
a representant of an element in
$ H_n^{f-{\rm proper}} (X (\RR), \RR)$,
such that  $f (\Gamma) $ is a bounded subset of $\RR_+$.
Assume the following assumptions hold~:
\begin{enumerate}
\item[(1)] The cohomology groups $H^i (X (\CC), f^{\ast} \cL_{\chi})$
are zero
for $i \not= n$ and almost all characters $\chi$
of finite order of $I_0$.
\item[(2)] The Euler characteristic $\chi (X (\CC))$ is equal to
$(- 1)^n$.
\end{enumerate}
Then
$$
\int_{\Gamma} f^s \omega = h (s) c^s \prod_d \Gamma (ds)^{\alpha (d)},
$$
where $h(s)$ is a non zero rational function in $\CC (s)$ and
$\alpha (d)$ is defined by
$$
\prod_d (1- T^d)^{\alpha (d) (-1)^n}
=
\frac{Z_0 (T)}{Z_{\infty} (T)},
$$
with $Z_0 (T)$ and $Z_{\infty} (T)$
respectively the zeta function of the monodromy action
around $0$ and  $\infty$ on $H^{\cdot}_c (f^{-1} (t), \CC)$ for $t$ a
generic point of
$\CC^{\times}$: with notations of \ref{zeta},
$Z_0 (T) = Z ([R f_{!} (\CC)_{\bar \eta_0}])$
and
$Z_{\infty} (T) = Z ([R f_{!} (\CC)_{\bar \eta_{\infty}}])$.
\end{theorem}

\begin{proof}Indeed, this follows from Theorem 4.2.10 and the remark
following it in 
\cite{LS}.
\end{proof}

Now we can deduce the following statement from formula (\ref{mdo}),
Theorem \ref{lsa}, and the material in \S\kern .15em 6.

\begin{theorem}\label{otherform}
The equality
$$
[(R \Delta_{\vert B !} \CC)^{G}_{\bar \eta_{0}}] -
[(R \Delta_{\vert B !} \CC)^{G}_{\bar \eta_{\infty}}] -
[(Rq^N_{|U !} \CC)^G_{\bar \eta_0}]
= (-1)^n
\sum_{i = 1}^n (V_{\phi^{d_i}} \otimes V_{d_i} - V_{\phi})
$$
holds in $K_{I_0}$.
\end{theorem}

\begin{proof}We  first write
\begin{equation}\label{newform}
\begin{split}
I (s + \frac{1}{2}) &=
\Gamma \Bigl(Ns + \frac{N + n}{2}\Bigr)
\int_{S} \Delta (x)^{s + \frac{1}{2}} 
\frac{\vert dx\vert}{\vert dq\vert}\\
&= |G|\Gamma \Bigl(Ns + \frac{N + n}{2}\Bigr)
\int_{\pi (S)} \Delta (y)^{s} 
\frac{\vert dy\vert}{\vert dy_1\vert},
\end{split}
\end{equation}
considering  the morphism $\pi~: \CC^n \rightarrow \CC^n / G = \CC^n$
sending $x$ to $y$ with $y_1 = q (x)$.
By the proof of Theorem \ref{exactbar}
and by 
Proposition \ref{bc}, we may  apply
Theorem \ref{lsa}
to $X = B / G$, $\Gamma = \pi (S)$ and $f$ the morphism induced by
$\Delta$.
Comparing with
(\ref{mdo}), one deduces
the relation (remind the duplication formula
$\Gamma (x + \frac{1}{2}) =
\sqrt{\pi} 2^{- 2x + 1} \frac{\Gamma
(2x)}{\Gamma (x)}$)
\begin{equation}
\frac{Z_0 (T)}{Z_{\infty} (T)}
=
(1 - (-1)^{\frac{N + n}{2}}T^N)^{(-1)^n}
\Bigl(\prod_{i = 1}^n \frac{1 - (-T)^{d_i}}{1 - (-T)}\Bigr)^{(-1)^{n - 1}}.
\end{equation}
The result follows by
rewriting everything in terms of virtual representations of $I_0$,
using
the fact that, by 
Theorem \ref{exactbar}
and Proposition \ref{ab} (3),
$(1 - (-1)^{\frac{N + n}{2}}T^N)^{(-1)^n}$ is the 
zeta function of the monodromy action
around $0$ on the cohomology of the Milnor fiber of $q^N_{|U}$.
\end{proof}

\subsection{Conclusion}Using Theorem \ref{conn} and Theorem 
\ref{compl}, one sees that Theorem \ref{deg} and
Theorem \ref{otherform} are in fact equivalent. Hence our proof of
Theorem \ref{otherform}
yields a new proof of Theorem \ref{deg} and shows that
Theorem \ref{deg} is in fact equivalent to knowing the precise form of
the gamma factors in Macdonald's formula.

\subsection*{}
{\small For the convenience of the reader, let us indicate for
each result which has not been proved in \cite{Trans}
the precise place in the present paper
where it is proved:
Assertions (1), (2), (3) and (4) in Proposition 3.2.1 of \cite{Trans}
are proved respectively in Proposition \ref{three}, Theorem \ref{compl},
Theorem \ref{conn} and Proposition \ref{ab}. 
Proposition 3.2.3 of \cite{Trans} is proved in Theorem \ref{exactbar},
Theorem 3.3 of \cite{Trans} in Remark \ref{last},
formula (0.5) of \cite{Trans} in Theorem \ref{max},
the assertion three lines above Lemma 4.3.3 of
\cite{Trans} in Proposition \ref{bc}.
}
\bibliographystyle{amsplain}

\end{document}